\documentclass[11pt]{article}

\begin{document}
 
\begin{center}

 {\bf Kolmogorov and Aleksandrov in Sevan Monastery, Armenia, 1929}

\vspace{0.1in}

V.G.Gurzadyan

\vspace{0.1in}

Yerevan Physics Institute, Armenia and University of Rome 'La Sapienza', Italy:
E-mail: gurzadya@icra.it

\vspace{0.1in}

(Published in {\it "The Mathematical Intelligencer}, vol.26, p.40-43, 2004.)

\vspace{0.2in}

\end{center}

Lake Sevan is one of beautiful spots in Armenia. It has always attracted poets and painters. Poet Maxim Gorky called it 'a piece of sky which came down among the mountains'. 
In the summer of 1929, Andrei N.Kolmogorov, 26, and Pavel S.Aleksandrov, 33, two mathematicians just starting their life-long friendship, visited Sevan. During about a month they enjoyed the beauties of the lake, living in a cell in the monastery on the island of Sevan. They visited other sites of Armenia, but actively continued their research. Both mathematicians mentioned Sevan in their memoirs, many decades later: 

\vspace{0.1in}

{\it At that time, Sevan was in its full glory. We bathed several times a day in its cool limpid waters. We also did a lot of work; in particular, I worked on "Aleksandrov-Hopf" [3] (typewriter was always with me) (Aleksandrov [1]).

\vspace{0.1in}

Naturally we were immediately attracted by the rocky islet, which with the drop in the level of the lake Sevan has now became a peninsula, and we wanted to settle there. This proved to be a simple matter. The cells of the monastery were empty and we occupied one of them (Kolmogorov [2]).} 

\vspace{0.1in}

Sevan provides an unusual combination of natural conditions, a high-altitude lake, 1900 m above sea level, 75 km in length and 55 km in breadth. The water is fresh and transparent, and the people of the region drink it abundantly. A number of rivers flow into the lake and only one, the Hrazdan, flows from it.  The lake is famous for its trout and beautiful peninsula which was once an island. During the Middle Ages various battles occurred over that rocky island, where Armenians often sheltered from their enemies.

The monastery on the island of Sevan was founded in 874 AD by the Armenian king Ashot, of the Bagratuni dynasty and princess Mariam. It includes two churches, the Arakelotz (the Church of Apostles) and Astvatsatsin (the Church of the Virgin) [4]. There is a panoramic view on the lake from the monastery.
In the 1930s, the water flow from the lake via the river Hrazdan was increased for irrigation and hydro-electric power purposes. This lowered the level of the lake, and transformed the island to a peninsula.  However, in 1929 the island still existed, and Kolmogorov colorfully describes the life on that isolated spot:
\vspace{0.1in}  

{\it The permanent population of the island consisted then of the archimandrite of the monastery (who had a fairly big house), his wife (who looked after some cows), the head of the meteorological station with his small family, and finally of the 'Captain' who did indeed command "the Sevan fleet", consisting of one motor boat and a few of the unusual pleasure boats. His picturesque figure has occasionally been described in literature (for example by Marietta Shaginyan).

Every day the archimandrite opened the lower church (there were two abandoned temples on the top of the hill), lit candles and, in complete solitude, recited the service.   Obviously the head of the meteorological station carried out his duties.  At times, the Captain brought honored guests, for example, Sar'yan or the then President of the All-Union Central Executive Committee of Armenia, but he was also ready to support humble tourists.

On the island, we both set to work.   With our manuscripts, typewriter, and a folding table, we sought out the secluded bays.   In the intervals between our studies, we bathed a lot.  To study I took refuge in the shade, while Aleksandrov lay for hours in full sunlight wearing only dark glasses and a white panama.  He kept this habit of working completely naked under the burning sun well into his old age.
On Sevan Aleksandrov worked on various chapters of his joint monograph with Hopf, "Topology" [3],  and he helped me to write the German text of my article on the theory of integrals.   Besides writing this paper, I was busy with ideas about the analytic description of Markov processes with continuous time, the end product of which later became the memoir "On analytic methods in the theory of probability".

Given its position Lake Sevan mostly enjoyed sunny weather, but sometimes clouds coming from the East filled up from the mountains, dropped down to the water and then, on contact with it, vanished.   We stayed there for about 20 days without leaving the place (apart from excursions to the monastery of Hayravank under the guidance of our captain).} 

\vspace{0.1in}

Hayravank monastery with its IX-century church is situated on a high rock facing Lake Sevan.   According his student, Grant Maranjian, Kolmogorov recalled tiny details on that stay at Sevan's island during his visit in Armenia in 1973. 

\vspace{0.1in}

{\it For Aleksandrov the day was approaching on which he had arranged to journey to Gagra, and we set off together for Yerevan (where we stayed for some days in a student hostel).  The temperature was 40°C, the sky was a hazy blue, and only after sunset did there unexpectedly appear the peak of Ararat, suspended in this blue sky.   We visited Echmiadzin (where we decided not to visit Katholikos as we did not have the right clothes).   From Echmiadzin we walked to Alagez (spending the first night by the lake, where the physicists who were working on cosmic showers very kindly put us up). After spending one night there (still without suits, wearing only shorts), we climbed the south summit of the Alagez, which did not present any complications (4000m).  From the top there opened up a view of the rocky northern summit (4100m), separated from the south by a huge ridge of snow, at the very bottom of which could be seen a small lake, its shores frozen and covered with snow.  Of course, Aleksandrov wanted to climb down there and bathe, but I preferred climbing the northern summit.}

\vspace{0.1in}

At the time, reaching the top of Mt.Aragatz (Alagez) involved a 30 km walk over mountainous terrain, the roads not having been built until the 1940s. Aragatz has a spectacular 3-km-wide crater with a glacier and 4 summits surrounding it. The northern summit is not only the highest but also the most difficult one for climbing; its easiest route is classified "1B" mountaineering category of difficulty.  

The cosmic ray station of Yerevan Physics Institute is situated at 3200 m above sea level. The nearby beautiful lake, mentioned by Kolmogorov, is an artificial reservoir that was constructed in the Second Millennium BC, as archaeological studies have found [5].  The surface of the lake is free of ice only during 2 to 3 months a year, and Aleksandrov's desire to bathe in the icy water is as characteristic, as  Kolmogorov's climbing of the most difficult northern summit of Aragatz.

After the dissolution of the Soviet Union, the former communist constraints on religion have been removed in the Republic of Armenia. The activity of the Sevan monastery is much increased, not only because of the larger number of priests but also because a seminary has been opened on the peninsula.
 
The monastery and the peninsula are now attractive touristic areas. Few of the tourists, however, can guess that the gloomy cells, with walls completely blackened by many centuries of candle smoke, once accommodated two leading mathematicians of the twentieth century.

\vspace{0.1in}

REFERENCES

\vspace{0.1in}

[1] P.S.Aleksandrov, {\it Russian Math.Surveys}, 35, 315, 1980. 

[2] A.N.Kolmogorov, {\it Russian Math.Surveys}, 41, 225, 1986.

[3] P.Alexandroff, H.Hopf, {\it Topologie}, Berlin, 1935. 

[4] P.Cuneo, {\it Architettura Armena}, De Luca Edittore, Rome, 1988.

[5] A.Kalantar, {\it Armenia: From the Stone Age to the Middle Ages}, Civilisations du Proche-Orient, Neuchatel-Paris, 1994.

\vspace{0.2in}

{\it Photo captions.}

1. Pavel Aleksandrov and Andrei Kolmogorov (photo Uspekhi Mat.Nauk, 1986).

2. The Sevan monastery where Kolmogorov and Aleksandrov lived in 1929 (picture by V.G.Gurzadyan, August 2003).

3. General view of the Sevan peninsula, in 1929 still an island, where the monastery is situated (photo by H.Badalian).

4. The summit of  Aragatz, climbed by Kolmogorov (photo H.Badalian).

5. The cosmic ray station, 3200m, where Kolmogorov and Aleksandrov stayed the night before climbing Aragatz mountain the next day. Many scientists who went mountaineering later visited the station.  The photo shows:  (from right to left), 
L.D.Faddeev (St.Petersburg), Alexander Migdal (Princeton), R.Kallosh (Stanford), Arcady Migdal (Landau Institute, Moscow), A.M.Polyakov (Princeton), V.G.Gurzadyan (Yerevan Physics Institute), A.G. Sedrakian (Yerevan Physics Institute) descending from the Aragatz summits in 1983.
Walking behind was another member of the group, A.Linde, now a famous cosmologist (his hat is visible in the back) (photo courtesy of V.G.Gurzadyan).

\end{document}